\documentclass[11pt]{amsart}
\newtheorem{theorem}{Theorem}

\newtheorem{lemma}[theorem]{Lemma}
\newtheorem{corollary}[theorem]{Corollary}
\theoremstyle{definition}

\begin{document}

\title[The logarithmic Sobolev inequality]{The logarithmic Sobolev inequality for a submanifold in Euclidean space}
\author{Simon Brendle}
\address{Department of Mathematics \\ Columbia University \\ New York NY 10027}
\begin{abstract}
We prove a sharp logarithmic Sobolev inequality which holds for submanifolds in Euclidean space of arbitrary dimension and codimension. Like the Michael-Simon Sobolev inequality, this inequality includes a term involving the mean curvature. 
\end{abstract}
\thanks{This project was supported by the National Science Foundation under grant DMS-1806190 and by the Simons Foundation.}

\maketitle

\section{Introduction}

In the early 1970s, Allard and, independently, Michael and Simon proved a Sobolev inequality which holds on every submanifold in Euclidean space (see \cite{Allard}, Section 7, and \cite{Michael-Simon}). This inequality is particularly useful on a minimal submanifold; in general, it contains a term involving the mean curvature. The constant in the Michael-Simon Sobolev inequality depends only on the dimension; however, the constant is not sharp. Castillon \cite{Castillon} gave an alternative proof of this inequality based on optimal transport. In a recent paper \cite{Brendle}, we proved a sharp version of the Michael-Simon Sobolev inequality for submanifolds of codimension at most $2$. In particular, this implies a sharp isoperimetric inequality for minimal submanifolds in Euclidean space of codimension at most $2$. 

In 2000, Ecker \cite{Ecker} proved a logarithmic Sobolev inequality which holds on every submanifold in Euclidean space. The constant in Ecker's inequality depends only on the dimension, but is not sharp. In this paper, we improve this inequality, and obtain the optimal constant. The following is the main result of this paper:

\begin{theorem} 
\label{main.thm.v1}
Let $\Sigma$ be a compact $n$-dimensional submanifold of $\mathbb{R}^{n+m}$ without boundary. Let $f$ be a positive smooth function on $\Sigma$. Then
\begin{align*} 
&\int_\Sigma f \, \Big ( \log f + n + \frac{n}{2} \, \log(4\pi) \Big ) \, d\text{\rm vol} - \int_\Sigma \frac{|\nabla^\Sigma f|^2}{f} \, d\text{\rm vol} - \int_\Sigma f \, |H|^2 \, d\text{\rm vol} \\ 
&\leq \bigg ( \int_\Sigma f \, d\text{\rm vol} \bigg ) \, \log \bigg ( \int_\Sigma f \, d\text{\rm vol} \bigg ), 
\end{align*} 
where $H$ denotes the mean curvature vector of $\Sigma$.
\end{theorem}

If we write $f = (4\pi)^{-\frac{n}{2}} \, e^{-\frac{|x|^2}{4}} \, \varphi$, then Theorem \ref{main.thm.v1} takes the following equivalent form:

\begin{corollary} 
\label{main.thm.v2}
Let $\Sigma$ be a compact $n$-dimensional submanifold of $\mathbb{R}^{n+m}$ without boundary, and let 
\[d\gamma = (4\pi)^{-\frac{n}{2}} \, e^{-\frac{|x|^2}{4}} \, d\text{\rm vol}\] 
denote the Gaussian measure on $\Sigma$. Let $\varphi$ be a positive smooth function on $\Sigma$. Then
\begin{align*} 
&\int_\Sigma \varphi \, \log \varphi \, d\gamma - \int_\Sigma \frac{|\nabla^\Sigma \varphi|^2}{\varphi} \, d\gamma - \int_\Sigma \varphi \, \Big | H + \frac{x^\perp}{2} \Big |^2 \, d\gamma \\ 
&\leq \bigg ( \int_\Sigma \varphi \, d\gamma \bigg ) \, \log \bigg ( \int_\Sigma \varphi \, d\gamma \bigg ), 
\end{align*} 
where $H$ denotes the mean curvature vector of $\Sigma$.
\end{corollary}

The logarithmic Sobolev inequality in Euclidean space has been studied by numerous authors; see e.g. \cite{Bakry}, \cite{Bobkov-Ledoux}, \cite{Gross}, \cite{Ledoux}. The Euclidean logarithmic Sobolev inequality can be proved in many different ways; for example, it can be deduced from the isoperimetric inequality in Gauss space (see \cite{Borell1}, \cite{Borell2}, \cite{Sudakov-Tsirelson}). The logarithmic Sobolev inequality also plays a central role in Perelman's monotonicity formula for the entropy under Ricci flow (see \cite{Perelman}).

The Gaussian measure $\gamma$ is very natural from a geometric point of view. In particular, the quantity $\gamma(\Sigma) = (4\pi)^{-\frac{n}{2}} \int_\Sigma e^{-\frac{|x|^2}{4}} \, d\text{\rm vol}$ is the Gaussian area which appears in Huisken's monotonicity formula for mean curvature flow (cf. \cite{Huisken}). The term $H + \frac{x^\perp}{2}$ vanishes if and only if $\Sigma$ is a critical point of the Gaussian area $\gamma(\Sigma)$. Submanifolds with this property are referred to as self-similar shrinkers: they shrink homothetically under mean curvature flow. We note that the Gaussian area and Huisken's monotonicity formula are of fundamental importance in the study of mean curvature flow; see e.g. \cite{Colding-Minicozzi}, \cite{Huisken}. 

The proof of Theorem \ref{main.thm.v1} follows a similar strategy as in \cite{Brendle} and is inspired in part by the Alexandrov-Bakelman-Pucci maximum principle (cf. \cite{Cabre}, \cite{Trudinger}). Alternatively, we could employ ideas from optimal transport; in that case, we would consider the transport map from $\mathbb{R}^{n+m}$ equipped with the Gaussian measure to the submanifold $\Sigma$ equipped with the measure $f \, d\text{\rm vol} = \varphi \, d\gamma$.

\section{Proof of Theorem \ref{main.thm.v1}}

Let $\Sigma$ be a compact $n$-dimensional submanifold of $\mathbb{R}^{n+m}$ without boundary. We denote by $T_x \Sigma$ the tangent space to $\Sigma$ at a point $x \in \Sigma$. Moreover, $T_x^\perp \Sigma$ will denote the normal space at $x \in \Sigma$. We denote by $I\!I$ the second fundamental form of $\Sigma$. In other words, given two tangent vector fields $X,Y$ and a normal vector field $V$, we define $\langle I\!I(X,Y),V \rangle = \langle \bar{D}_X Y,V \rangle = -\langle \bar{D}_X V,Y \rangle$, where $\bar{D}$ denotes the connection on $\mathbb{R}^{n+m}$. Finally, the mean curvature vector $H$ is defined as the trace of the second fundamental form $I\!I$.

We now give the proof of Theorem \ref{main.thm.v1}. We first consider the special case that $\Sigma$ is connected. By scaling, we may assume that 
\[\int_\Sigma f \, \log f \, d\text{\rm vol} - \int_\Sigma \frac{|\nabla^\Sigma f|^2}{f} \, d\text{\rm vol} - \int_\Sigma f \, |H|^2 \, d\text{\rm vol} = 0.\] 
Since $\Sigma$ is connected, the operator $u \mapsto \text{\rm div}_\Sigma(f \, \nabla^\Sigma u)$ has one-dimensional kernel and cokernel. Hence, we can find a smooth function $u: \Sigma \to \mathbb{R}$ such that 
\[\text{\rm div}_\Sigma(f \, \nabla^\Sigma u) = f \, \log f - \frac{|\nabla^\Sigma f|^2}{f} - f \, |H|^2.\] 
We define 
\begin{align*} 
U &:= \{(x,y): x \in \Sigma, \, y \in T_x^\perp \Sigma\}, \\ 
A &:= \{(x,y) \in U: D_\Sigma^2 u(x) - \langle I\!I(x),y \rangle \geq 0\}. 
\end{align*} 
Moreover, we define a map $\Phi: U \to \mathbb{R}^{n+m}$ by 
\[\Phi(x,y) = \nabla^\Sigma u(x) + y\] 
for all $(x,y) \in U$. Note that $|\Phi(x,y)|^2 = |\nabla^\Sigma u(x)|^2+|y|^2$ since $\nabla^\Sigma u(x) \in T_x \Sigma$ is orthogonal to $y \in T_x^\perp \Sigma$. \\

\begin{lemma}
\label{Phi.surjective} 
We have $\Phi(A) = \mathbb{R}^{n+m}$. 
\end{lemma}

\textbf{Proof.} 
Given $\xi \in \mathbb{R}^{n+m}$, we can find a point $\bar{x} \in \Sigma$ where the function $w(x) := u(x) - \langle x,\xi \rangle$ attains its minimum. Clearly, $\nabla^\Sigma w(\bar{x}) = 0$. From this, we deduce that $\xi = \nabla^\Sigma u(\bar{x}) + \bar{y}$ for some vector $\bar{y} \in T_{\bar{x}}^\perp \Sigma$. Moreover, $D_\Sigma^2 w(\bar{x}) \geq 0$. This implies $D_\Sigma^2 u(\bar{x}) - \langle I\!I(\bar{x}),\xi \rangle \geq 0$, and consequently $D_\Sigma^2 u(\bar{x}) - \langle I\!I(\bar{x}),\bar{y} \rangle \geq 0$. Thus, $(\bar{x},\bar{y}) \in A$ and $\Phi(\bar{x},\bar{y}) = \xi$. \\

\begin{lemma}
\label{Jacobian.determinant}
The Jacobian determinant of $\Phi$ is given by 
\[\det D\Phi(x,y) = \det (D_\Sigma^2 u(x) - \langle I\!I(x),y \rangle)\] 
for each point $(x,y) \in U$.
\end{lemma}

\textbf{Proof.} 
See \cite{Brendle}, Lemma 5. \\

\begin{lemma}
\label{Jacobian.estimate}
The Jacobian determinant of $\Phi$ satisfies 
\[0 \leq e^{-\frac{|\Phi(x,y)|^2}{4}} \, \det D\Phi(x,y) \leq f(x) \, e^{-\frac{|2H(x)+y|^2}{4} - n}\] 
for each point $(x,y) \in A$.
\end{lemma}

\textbf{Proof.} 
Consider a point $(x,y) \in A$. Using the identity $\text{\rm div}_\Sigma(f \, \nabla^\Sigma u) = f \, \log f - \frac{|\nabla^\Sigma f|^2}{f} - f \, |H|^2$, we obtain 
\begin{align*} 
\Delta_\Sigma u(x) - \langle H(x),y \rangle 
&= \log f(x) - \frac{|\nabla^\Sigma f(x)|^2}{f(x)^2} - |H(x)|^2 \\ 
&- \frac{\langle \nabla^\Sigma f(x),\nabla^\Sigma u(x) \rangle}{f(x)} - \langle H(x),y \rangle \\ 
&= \log f(x) + \frac{|\nabla^\Sigma u(x)|^2+|y|^2}{4} \\ 
&- \frac{|2 \, \nabla^\Sigma f(x) + f(x) \, \nabla^\Sigma u(x)|^2}{4 f(x)^2} - \frac{|2H(x)+y|^2}{4} \\ 
&\leq \log f(x) + \frac{|\nabla^\Sigma u(x)|^2+|y|^2}{4} - \frac{|2 H(x)+y|^2}{4}.
\end{align*}
Moreover, $D_\Sigma^2 u(x) - \langle I\!I(x),y \rangle \geq 0$ since $(x,y) \in A$. Using the elementary inequality $\lambda \leq e^{\lambda-1}$, we obtain 
\begin{align*} 
0 &\leq \det (D_\Sigma^2 u(x) - \langle I\!I(x),y \rangle) \\ 
&\leq e^{\text{\rm tr}(D_\Sigma^2 u(x) - \langle I\!I(x),y \rangle) - n} \\ 
&\leq f(x) \, e^{\frac{|\nabla^\Sigma u(x)|^2+|y|^2}{4} - \frac{|2H(x)+y|^2}{4} - n}.
\end{align*}
Using Lemma \ref{Jacobian.determinant} and the identity $|\Phi(x,y)|^2 = |\nabla^\Sigma u(x)|^2+|y|^2$, we conclude that 
\[0 \leq \det D\Phi(x,y) \leq f(x) \, e^{\frac{|\Phi(x,y)|^2}{4} - \frac{|2H(x)+y|^2}{4} - n}.\] 
This proves the assertion. \\

After these preparations, we now complete the proof of Theorem \ref{main.thm.v1}. Using Lemma \ref{Phi.surjective} and Lemma \ref{Jacobian.estimate}, we obtain 
\begin{align*} 
1 &= (4\pi)^{-\frac{n+m}{2}} \int_{\mathbb{R}^{n+m}} e^{-\frac{|\xi|^2}{4}} \, d\xi \\ 
&\leq (4\pi)^{-\frac{n+m}{2}} \int_\Sigma \bigg ( \int_{T_x^\perp \Sigma} e^{-\frac{|\Phi(x,y)|^2}{4}} \, |\det D\Phi(x,y)| \, 1_A(x,y) \, dy \bigg ) \, d\text{\rm vol}(x) \\ 
&\leq (4\pi)^{-\frac{n+m}{2}} \int_\Sigma \bigg ( \int_{T_x^\perp \Sigma} f(x) \, e^{-\frac{|2H(x)+y|^2}{4} - n} \, dy \bigg ) \, d\text{\rm vol}(x) \\ 
&= (4\pi)^{-\frac{n}{2}} \, e^{-n} \int_\Sigma f(x) \, d\text{\rm vol}(x).
\end{align*}
Consequently, 
\[n + \frac{n}{2} \, \log(4\pi) \leq \log \bigg ( \int_\Sigma f \, d\text{\rm vol} \bigg ).\] 
Combining this inequality with the normalization 
\[\int_\Sigma f \, \log f \, d\text{\rm vol} - \int_\Sigma \frac{|\nabla^\Sigma f|^2}{f} \, d\text{\rm vol} - \int_\Sigma f \, |H|^2 \, d\text{\rm vol} = 0\] 
gives 
\begin{align*} 
&\int_\Sigma f \, \Big ( \log f + n + \frac{n}{2} \, \log(4\pi) \Big ) \, d\text{\rm vol} - \int_\Sigma \frac{|\nabla^\Sigma f|^2}{f} \, d\text{\rm vol} - \int_\Sigma f \, |H|^2 \, d\text{\rm vol} \\ 
&= \int_\Sigma f \, \Big ( n+\frac{n}{2} \, \log(4\pi) \Big ) \, d \text{\rm vol} \\ 
&\leq \bigg ( \int_\Sigma f \, d\text{\rm vol} \bigg ) \, \log \bigg ( \int_\Sigma f \, d\text{\rm vol} \bigg ). 
\end{align*} 
This proves Theorem \ref{main.thm.v1} in the special case when $\Sigma$ is connected.

It remains to consider the case when $\Sigma$ is disconnected. In that case, we apply the inequality to each individual connected component of $\Sigma$, and sum over all connected components. Since 
\[a \log a + b \log b < a \log(a+b) + b \log (a+b) = (a+b) \log (a+b)\] 
for $a,b>0$, we conclude that 
\begin{align*} 
&\int_\Sigma f \, \Big ( \log f + n + \frac{n}{2} \, \log(4\pi) \Big ) \, d\text{\rm vol} - \int_\Sigma \frac{|\nabla^\Sigma f|^2}{f} \, d\text{\rm vol} - \int_\Sigma f \, |H|^2 \, d\text{\rm vol} \\ 
&< \bigg ( \int_\Sigma f \, d\text{\rm vol} \bigg ) \, \log \bigg ( \int_\Sigma f \, d\text{\rm vol} \bigg ) 
\end{align*} 
if $\Sigma$ is disconnected. This completes the proof of Theorem \ref{main.thm.v1}.

\section{Proof of Corollary \ref{main.thm.v2}}

In this section, we explain how Corollary \ref{main.thm.v2} follows from Theorem \ref{main.thm.v1}. We consider the Gaussian measure $d\gamma = (4\pi)^{-\frac{n}{2}} \, e^{-\frac{|x|^2}{4}} \, d\text{\rm vol}$ on $\Sigma$. Moreover, we write $f = (4\pi)^{-\frac{n}{2}} \, e^{-\frac{|x|^2}{4}} \, \varphi$. Clearly, $f \, d\text{\rm vol} = \varphi \, d\gamma$. Using the identity $\frac{\nabla^\Sigma f}{f} = \frac{\nabla^\Sigma \varphi}{\varphi} - \frac{x^{\text{\rm tan}}}{2}$ and the relation $\text{\rm div}_\Sigma(x^{\text{\rm tan}}) = n + \langle H,x^\perp \rangle$, we obtain 
\begin{align*} 
&\log f + n + \frac{n}{2} \, \log(4\pi) - \frac{|\nabla^\Sigma f|^2}{f^2} - |H|^2 - \frac{\text{\rm div}_\Sigma(f \, x^{\text{\rm tan}})}{f} \\ 
&= \log \varphi - \frac{|\nabla^\Sigma \varphi|^2}{\varphi^2} - \Big | H + \frac{x^\perp}{2} \Big |^2. 
\end{align*} 
In the next step, we integrate the expression on the left hand side with respect to the measure $f \, d\text{\rm vol}$, and we integrate the expression on the right hand side with respect to the measure $\varphi \, d\gamma$. Using the divergence theorem, we conclude that 
\begin{align*} 
&\int_\Sigma f \, \Big ( \log f + n + \frac{n}{2} \, \log(4\pi) \Big ) \, d\text{\rm vol} - \int_\Sigma \frac{|\nabla^\Sigma f|^2}{f} \, d\text{\rm vol} - \int_\Sigma f \, |H|^2 \, d\text{\rm vol} \\ 
&= \int_\Sigma \varphi \, \log \varphi \, d\gamma - \int_\Sigma \frac{|\nabla^\Sigma \varphi|^2}{\varphi} \, d\gamma - \int_\Sigma \varphi \, \Big | H + \frac{x^\perp}{2} \Big |^2 \, d\gamma. 
\end{align*}
Hence, Corollary \ref{main.thm.v2} follows from Theorem \ref{main.thm.v1}.

\end{document}